\documentclass[11pt]{article}

\usepackage{epsfig}
\usepackage{graphicx}
\usepackage{color}

\newcommand{\DS}{\LARGE \renewcommand{\baselinestretch}{1.67} \normalsize}

\newtheorem{theorem}{Theorem}[section]

\title{Algorithms for Deforming and Contracting Simply Connected Discrete
Closed Manifolds (II) }

\author{
Li Chen \\
Department of Computer Science and Information Technology\\
University of the District of Columbia\\
Email: lchen@udc.edu\\
}

\begin{document}
\DS

\maketitle

\abstract
{In an exploration paper, {\it L. Chen, Algorithms for Deforming and Contracting Simply Connected Discrete
Closed Manifolds (I)},  we designed algorithms for
 deforming and contracting a simply connected discrete
closed manifold to a discrete sphere.  However, the algorithms could not guarantee to be applicable to every case.
This paper will be the continuation of the exploration.

This paper contains two main procedures: (1) A shrinking procedure to contract a simply connected closed manifold. Unlike ones in
the previous paper, we added a tree structure to support the process. (2) A more direct procedure for
mapping a component from a separated simply connected closed manifold to a disk.

We also discuss the practical use of these algorithms in topological data analysis. We think that we
have an algorithmic solution, but careful detailed analysis should be done next.
 
}

\section{Introduction}

We recall some basic results in this section. In general, any smooth real $m$-dimensional manifold can be smoothly embedded in $R^{2m}$; this is called the (strong) Whitney Embedding Theorem. And any $m$-manifold with a Riemannian metric (Riemannian manifold) can be isometrically imbedded to
an $n$-Euclidean space, where $n\le c\cdot m^3$, $c$ is small constant. This is called the Nash Embedding Theorem. Therefore, we can discuss our problem in Euclidean space or a space that can be easily embedded to Euclidean space.

On the other hand, according to  Whitehead \cite {Whit40}:  Every smooth manifold admits an (essentially unique) compatible piecewise linear structure.
In 1952, Moise proved the following theorem \cite{Moi77}:  Any 3-dimensional manifold is smooth, and thus piecewise linear.

Therefore, we can just discuss discrete manifolds in a partitioned Euclidean space for any type of smooth manifolds.

Our new method will be based on the previous paper, L. Chen, Algorithms for Deforming and Contracting Simply Connected Discrete
Closed Manifolds (I), {\it https://arxiv.org/abs/1507.07171}.

In this paper, we will do the following: (1)
For a closed and simply connected discrete $m$-manifold $M_m$ in $n$-dimensional Euclidean space $E_n=R^n$ , we will fill a discrete $(m+1)$-manifold that is
bounded by $M_m$. If this filling is valid, then we will design an algorithm that can contract the boundary of $F$ to be the boundary of a single
$(k+1)$-cell that is homeomorphic to a $k$-sphere.  This paper will add a tree structure in the algorithm.   
(2) A more direct procedure for
mapping a component from a separated simply connected closed manifold to a disk. 

This paper is still an exploration paper.

\section{Some Reviews}

In \cite{Che15I}, we observed that Chen-Krantz actually proved the following result: A simply connected (orientable) manifold $M$ in space $U$.  If $M$ is a supper submanifold, the dimension of $M$ is smaller than the dimension of $U$ by one, in such a case, we can use Jordan's theorem to first separate the $U$ into two components. The deformation becomes the pure contraction. This result can be obtained directly from Chen-Krantz's paper. (L. Chen and S. Krantz,  A Discrete Proof of The General Jordan-Schoenflies Theorem,
{\tt http://arxiv.org/abs/1504.05263})

However, if the dimension of $U$ is much bigger than
the dimension of $M$, we will need other ways, for instance, we need to fill an $(m+1)$-manifold bounded by $M$, where $m$ is the dimension of $M$. Some algorithms have been discussed in \cite{Che15I}. But these algorithms may not work for some cases.

In this paper, we continue the task  of finding the way of filling of $M$ and also discuss a method of deduction the cells on $M$.

\section{Two Algorithms for the Closed Simply Connected Manifolds}

In this section, we present two algorithms for the closed simply connected manifolds.

\subsection{The Filling Procedure for Simply Connected Manifolds}

In this section, we will continue our discussion in the previous paper (I) . We still want to find a $(m+1)$ dimensional filling
of $M$  in $E_n$ \cite{Che15I}.

Let $U=E_n$ be the $n$-dimensional Euclidean Space. $\Sigma_n$ is a PL decomposition of $E_n$. More specifically, $\Sigma_n$ can be a cubic, simplicial, or other discrete decomposition of $E_n$ discussed in \cite{Che14}.

$M_m$ be a simply connected discrete $m$-manifold in $\Sigma_n$. $M_m$ is closed and orientable.  we will need other ways, for instance, we need to fill an $(m+1)$-manifold bounded by $M_m$, where $m$ is the dimension of $M_m$.

We define a branch is a connected component of $M_m$ and the component will contain at least one point that has a local maximum positive curvature (positive sectional curvature for each dimension). We call such an area a peak.

The following algorithm will use a tree-structure to record a branch (Fig.1). And the total tree will represent the branch structures
of the discrete manifold.  The tree structure will provide algorithmic advantages in real time calculation for filling.

\begin{figure}[hbt]

\begin{center}

 \includegraphics [width=3in] {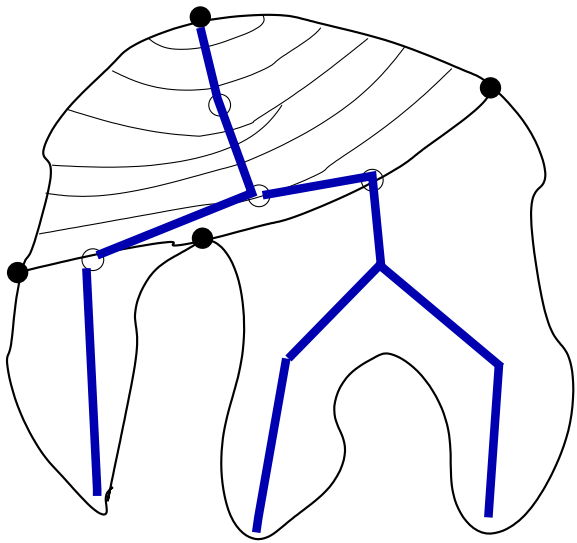} 

\end{center}
\caption{The general structure of the filling and its supporting tree. }
\end{figure}

{\bf Algorithm 3.A} . This algorithm is not the same as Algorithm 3.1 in (I)

\begin{description}
	
\item [Step 1] Make $M_m$ to be a local flat $m$-manifold in $\Sigma_n$. From the top (or left) direction of the minimum cubical box that contains $M_m$ in $\Sigma_n$. Find the first $m$-cell $e$ in $M_m$, remove (or mark) this cell. Obtain the boundary of $M_m-\{e\}$ ($M_m-\{e\}$ means to remove the inner part of $e$.) This boundary $B$ is an $(m-1)$-cycle. For instance, if $M_m$ is a closed curve, $B$ contains two points.

\item [Step 2]  There is an $(m+1)$-cell $g$ containing $e$ in $\Sigma_n$ ,  Boundary of g, $\partial g$ is an $m$-cycle. ${E} =\partial g-e$ is a collection of $m$-cells that have the same edge as $e$, $B$. We choose ${E}$ has the minimum $m$-cells. We know that ${\bar E}$ and $e$ are gradually varied.

\item [Step 3]  Find a $B'$ that is gradually varied to $B$ on $M_m$ with the largest different number of elements from $B$ ($B\ne B'$). Make a simply connected $m$-manifold with boundary $B'$, $E'$. Here is the priority: first we want  $E'$ and $E$ are gradually varied. If we can get one, we are looking for such a $E'$ has the smallest number of $m$-cells. If we could not get one t hat is

    If $E'$ contains the number of $m$-cells that is bigger than one that contains a set that is an old $B$ or can be a boundary that is the removed set or marked set.
    then, this $E'$ is not necessary. It means that we find a branch that is all removed or marked elements.

    If $E'$ must contain a isolated element or elements in $M_m$ , that means its time to make a branch(s).

\item [Step 4]  If $E'$ must contain a isolated element or elements in $M_m$ (Fig. 2 and Fig. 3.), there must be two cases.
First,  make a branch(s) if $E'\cap M_m$ is a connected two or more cycles, it can be cut one out (Fig.3). Second, if it is the union of two disconnected cycles, we can use the Step 1 - Step 3 to fill a branch(s) ( we define it as in the outer part) then to remove it; this equivalent to
push back the branch, this is just like a negative curvature point (Fig.2).

\item [Step 5] Now the only problem is to deal with the case that we cannot get two gradually varied $E$ and $E'$ . Using modular 2 sum of $E$ and $E'$ we can have a closed m-cycle(s) . This cycle(s) are much smaller. we can fill this cycle with the method from Step 1-Step 4 by inserting the gradually varied sub "$E$."

\item [Step 6]  The structure of the tree is based on the center of $E$ in each filling. When a branch is made, the cut (the last $E'$ in a subsequence) will have a potential link to the parent part. The node will be attached to the center of the cut.
    For complex case, we can use two trees, one is the inner tree and another is the outer tree (or set of outer trees called
    outer forest). Combining all together, we can make an (m+1)-manifold that has the boundary $M_m$.  (As we discussed in (I), if $n> m+1$ there might be multiple choices when actual perform this algorithm. )

\end{description}

See some situations shown in Fig. 2 and Fig. 3.  The $E$ is a minimum cap when removing cells on $M_m$. This process will determine a set of
gradually varied fills for each branch. After all, a (m+1)-manifold with the boundary that is $M_m$ will be determined.
Then we use the algorithm in (I) will be able to do a contraction . There are still some details need to be done in this algorithm.

\begin{figure}[hbt]

\begin{center}

 \includegraphics [width=3in] {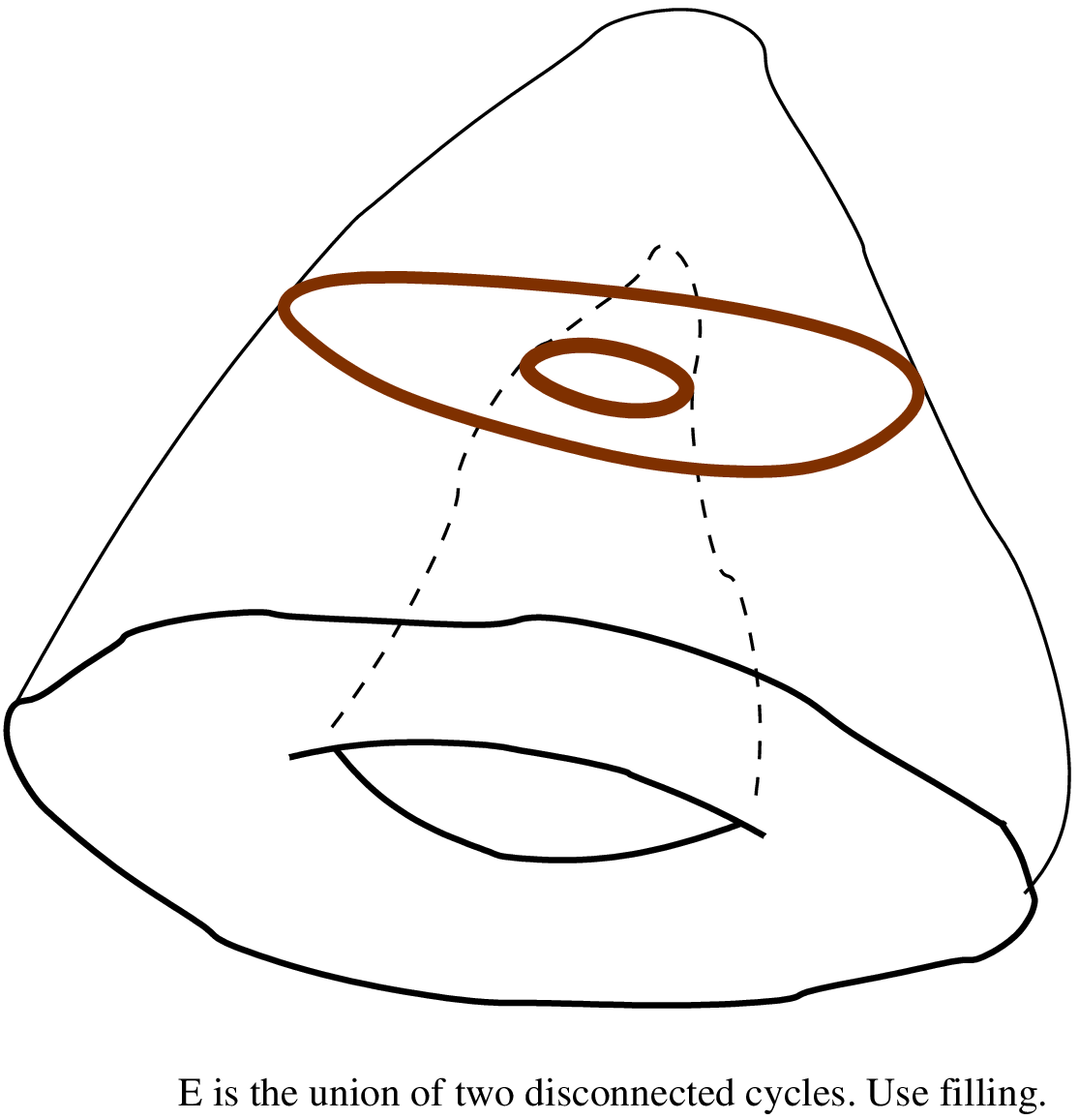} 

\end{center}
\caption{A case needs to fill other side of the surface. It is equivalent to push it done at the elliptical area.    }

\end{figure}

\begin{figure}[hbt]

\begin{center}

 \includegraphics [width=3in] {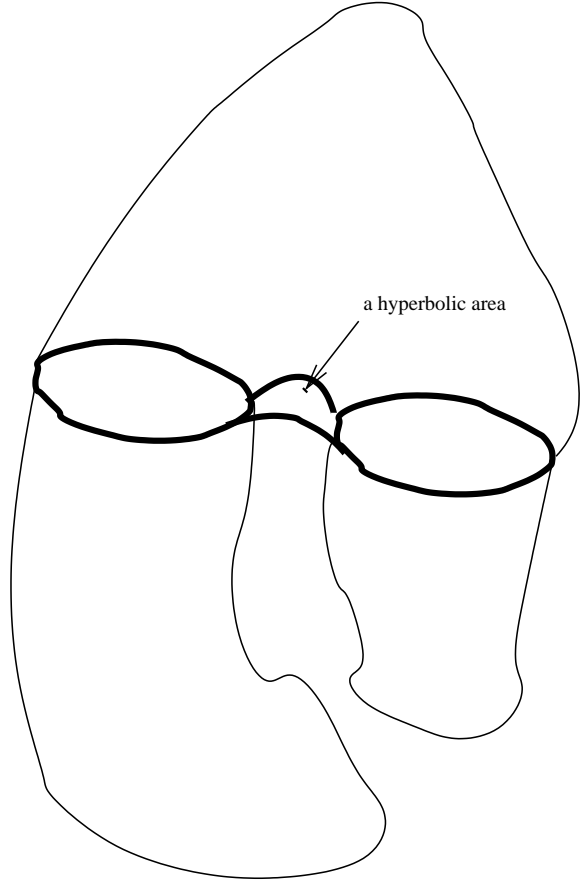} 

\end{center}
\caption{ A case needs to make a branch. }

\end{figure}

\subsection{The Reduction Procedure for Simple Connectedness }

Using the cell distance from A, and A' cell on M to choose the closer one m+1 cell D (

Based on a theorem (Theorem 5.1) proved  by Chen and Krantz ( L. Chen and S. Krantz,  A Discrete Proof of The General Jordan-Schoenflies Theorem, \newline {\tt http://arxiv.org/abs/1504.05263}), we concluded that a simply connected closed $m-1$-manifold split a simply connected closed $m$-manifold into two components, each of which is simply connected. We require both $B$ and $M$ are locally flat. We assume all manifolds discussed in this paper will be orientable.

We restate this theorem as follows :
\begin{theorem}
(The Jordan Theorem for the closed surface on a 3D manifold) Let $M_3$
be a simply connected 3D manifold (discrete or $PL$); a closed discrete surface $S$ (with
local flatness) will separate $M_3$ into two components. Here $M_3$ can be closed.
\end{theorem}

Based on this theorem we will design a procedure that will generate a homeomorphic mapping for a component in $M_3-S$ to a 3-disk.
So if $M_3$ is closed then $M_3$ is homeomorphic to a 3-sphere.)

\begin{theorem}
If $M_3$ is closed in Theorem 3.1, we can algorithmically make $M_3$ to be homeomorphic to a 3-sphere in discrete case.
\end{theorem}

{\bf The Algorithmic Proof:}


According to Theorem 3.1, we already proved that in discrete case,  a simply connected closed $2$-manifold (orientable) $B$ with local flatness will split a  $3$-manifold $M_3$ into two components $D$ and $D'$. We now show that each of them will be simply connected. 

In fact, each of the two components will be simply connected. This is because that if $M_3$ is simply connected. A simply closed curve $C$ is contractible to a point $p\in C$ on $M_3$.  Let $\Omega$ be the contraction sequence in discrete case, we call gradually variation in \cite{Che14}. We might as well assume $p$ is not on $B$.

The contraction sequence $\Omega$ may contain some point on $B$, we can modify the contraction by using $\Omega\cap B$ to replace
$\Omega$ to get a new contraction sequence. So this theorem is valid. A curve in $\Omega$ may intersect with $C$ but will not cross-over $C$.

A curve started in a component $D$ will pass (have both enter to $B$ and out of $B$ to another component called a pass)  even times on $B$ and also finite number of times in discrete case. When $C'\in \Omega$ pass $B$ and it will enter $B$ from $D'$ , we can find
a curve $C'_B(0)$ on $B$ two link two points (the last point of leaving of $B$ and the fist point entering $B$, since $B$ is simply connected.) Use $C'_B(0) $ to replace that corresponding arc in $D'$ . If there are multiple passes, we can use  Use $C'_B(1),...,C'_B(i) $ to replace all. So we get a new $C'_B$ that only contains points in $D\cup B$. Use the same
process for all curves in  $\Omega$, we will get a  $\Omega_B$. That is the contraction to $p$. So $D$ is simply connected and, so is $D'$

Note that If the curve only stay on $B$ and back to the original component $D$ will not be counted as a pass.

Next, we would design an algorithm to show that $D$ is homeomorphic to an $m$-disk.

Again,  according to the definition of simple connectedness, a closed 1-cycle will be contractible on $M$. since we assume that
there is no edge in $M$ , $M$ is closed and orientable, then there is no holes in $M$. (if there is a hole, then the edge of a hole will be a 1-cycle and it is not contractible. ). For instance, if $M$ is a torus, some cycles are not contractible.


Let $m=3$ for now, we will see $m$ can be any number.

Step (1): Remove an $m$-cell $e$ from $M_m$ will leave a $(m-1)$-cycle. This cycle $B$ is always simply connected as well. we use this property (plus the theorem we discussed above). (Note:  We can assume that any closed simply connected $(m-1)$-manifold is an $(m-1)$-cycle that is homeomorphic to $(m-1)$-sphere when we prove for $m$.) Let $D=M_m-e$ . Note that this subtraction is
to remove the $m$-cell $e$ not its $(m-1)$ edges (faces).

Step (2): We will remove more $m$-cells of $M_m$ if they have an edge(face) on $B$. Algorithmically, we remove another $e\in D$ and $e$ has an $(m-1)$-edge(face) in $B$. We also denote the new edge of $D$ to be $B$ as well.
$B$ will be a new edge set of $D=M-\{\mbox{removed $m$-cells}\}$.
$B$ is still the $(m-1)$-cycle (pseudo-manifold in \cite{Che14}) since use the boundary of new $e$ in $D-B$ to replace $e\cap B$ .

Step (2'): For actual design of the algorithm, we will calculate the cell-distances to all points in $M_m$ from a fix point $o$, $e\ne o$, determine beforehand. We always section new $e$ that is adjacent to $B$ has the greatest distance to $o$. This is a strategy for the balanced selection of new $e$.

Step (3): Since $B$ is an $(m-1)$-cycle, so $D$ is always simply connected based on Theorem 3.1.   $M_m$ only contains finite number of $m$-cells, It mean that this process will end at the $Star(o)$ in $M_m$. Therefore, we algorithmically showed that
in discrete space, $M_m-e$ is continuously shrinking (homomorphic) to $Star(o)\cap M_m$. Since $Star(o)\cap M_m$ is
an $m$-disk. So, the reversed steps determine a homomorphic mapping  from an $m$-disk to $M_m-e$.

We also know that $e$ is also homomorphic to $m$-disk. The connected-sum of $e$ and $M_m-e$ is homomorphic to $m$-sphere
where $m=3$. This is the end of the algorithmic {\bf proof}.

For some specially cases when local flatness are considered, we can use some sub-procedure to modify. Just like
we did it in the Chen-Krantz paper.

Using this result, we can design the algorithms or procedure to show that in discrete space, a closed simply connected $m$-manifold is homeomorphic to $m$-sphere.

\begin{theorem}
(The Jordan Theorem of Discrete $m$-manifolds) Let $M_m$
be a simply connected discrete or PL $m$-manifold; a closed simply connected discrete $(m-1)$-manifold $B$ (with
local flatness) will separate $M_m$ into two components. Here $M_m$ can be closed.
\end{theorem}

(With the recursive assumption, $B$ can be assumed to be  homeomorphic to an $(m-1)$-sphere.)

\begin{theorem}
 If $M_m$ is closed in Theorem 3.3, we can algorithmically make $M_m$ to be homeomorphic to a $m$-sphere in discrete case.
\end{theorem}

Some discussions:  In Step (2),  if $B$ is not a simple cycle, then $M$ is not simply connected according to the
theorem above. We will use this property to determine if $M$ is simply connected in the next section.

In Step (3), since $M$ contains finite number of $m$-cells, there will be definite always to reach the end.
It is also possible to use gradually varied "curve" or $(m-1)$-cycle on $M_m$ of $B$ to replace $B$ to reach the maximum number of removal in practice.  However, we have to keep the removal balanced to a certain point meaning that we try always to remove one that has the furthermost distance on the edge to the fixed point $o$ (which we contract to).

This entire process of the algorithm determined a homeomorphism to the $m$-sphere. Therefore, we would like to say that
in discrete space, a closed-orientable simply connected $m$-manifold is homomorphic to an $m$-sphere.

As an equivalent statement, we observed that a closed-orientable $m$-simply connected manifold $M_m$ is a homeomorphic to $m$-sphere if and only if there is $(m+1)$-disk that is simply connected and has the boundary that is $M_m$. The discussion is in the last subsection 3.2.

\section{The Algorithm for Determining Simple Connectedness}

This is a revised procedure. To decide  if a discrete $m$-manifold is simply connected, we can
use the procedure described in Section 3 .  If there is boundary $B$ that is a union of two or more $(m-1)$-cycles
That means $M_m$ is not a simply connected manifold.

The algorithm for deciding that a complex is an $m$-manifold was described in \cite{Che12,Che14}. So our algorithm
will be first decide if $M_m$ is an $m$-manifold.

{\bf Algorithm 4.A} The algorithm of deciding if a discrete $m$-manifold is simply connected. 
In the proof in Section 3.2, we already suggest
such a procedure . Here we only need to rewrite it as an algorithm. The key part of the algorithm is
to check every deleting of an $m$-cell in the procedure will maintain the boundary to be a single $(m-1)$-cycle
(a simply connected closed $(m-1)$-manifold).
In this algorithm, we assume that we have a set of all $i$-cells for the complex $M_m$, $0\le i \le m$. 
We have already checked that $M_m$ is a closed $m$-manifold. 

\begin{description}

\item [Step 1:] Define a point $o$ in $D=M_m$ as the origin. Calculate all cell-distances from $o$ to all cells.

\item [Step 2:] Remove an $m$-cell $e\in D$ that is furthermost from $o$. It will leave a $(m-1)$-cycle, $B$. This cycle $B$ is always simply connected as well.  In this section, the $(m-1)$-cycle and $(m-1)$-simple cycle will be different. $(m-1)$-cycle is now the closed $(m-1)$-cell path where two adjacent cells share a $(m-2)$-cell.

\item [Step 3:] We will remove a set of $m$-cells of $D$ they are adjacent to $B$. Algorithmically, we remove a new $e\in D$ that is adjacent to $B$ and it is furthermost from $o$. The cell-distance will be used to determine this distance. If the new boundary $B$ of $D$ is a simple closed path. We continue this step. Otherwise, we report $M_m$ is not simply connected.

\item [Step 3':] The procedure to decide a path is a simple closed path: Check if a cell is used more than once if it is not at the beginning or end of the path.

\end{description}

In topological data analysis, it is common to ask if a data set is simply connected. This algorithm will work with
the algorithms of deciding if a simplicial complex or cellular complex is a discrete manifold \cite{Che14}. After it is done,
we can apply it to decide if this manifold is simply connected.

\section{Conclusion}

In this paper, we showed a new way to deal with the closed $m$-manifold. It also could have 
some real world applications in topological data analysis. 

In discrete cases, a closed $M$ with simple connectedness is the same as any closed $(m-1)$-cycle (with local flatness) can separate $M$ into two disconnected component. Each component is homeomorphic to an $m$-disk.

It seems like that we
can prove algorithmically that a simply connected closed 3D manifold in discrete case is homeomorphic to a 3-sphere.
However, due to the fact that some special cases might be exist, we might need to do some rechecks.  

In other words, theoretically, this algorithmic proof should be carefully checked including with some actual coding for real world problems.

The algorithmic proof might differs other proofs since we only can deal with finite number of cells in the manifolds.

\section{Appendix: Some Concepts in Manifolds and Discrete Manifolds }

The concepts of this paper are in \cite{Che14}. We also use some concepts from the following two papers:
L. Chen    A Concise Proof of Discrete Jordan Curve Theorem, {\tt http://arxiv.org/abs/1411.4621} and
L. Chen and S. Krantz,  A Discrete Proof of The General Jordan-Schoenflies Theorem,
\newline {\tt http://arxiv.org/abs/1504.05263}.

{\it
A discrete space is a graph $G$ having an associated structure. We always assume that $G$ is finite, meaning that $G$ contains only a
finite number of vertices.  Specifically, ${\cal C}_2$ is the set of all minimal cycles representing all possible 2-cells;
$U_2$ is a subset of  ${\cal C}_2$.   Inductively, ${\cal C}_3$  is the set of all minimal 2-cycles made by $U_2$.   $U_3$ is a subset of  ${\cal C}_3$ .
Therefore $\langle G,U_2,U_3,\cdots,U_k \rangle$ is a discrete space. We can see that a simplicial complex is a discrete space.  For computational purposes,  we want to require
that each element in $U_i$ can be embedded into a Hausdorff space or Euclidean space using a polynomial time algorithm （or an efficient constructive method). And such a mapping will
be a homeomorphism to an $i$-disk with the internal area of the cell corresponding to an $i$-ball that can be determined also in polynomial time.
Another thing we need to point out here
is that $u\cap v$ in   $\langle G,U_2,U_3,\cdots,U_k \rangle$ must be connected. In most cases,  $u\cap v$ is a single $i$-cell in $U_i$ or empty.
In general,  $u\cap v$ is homeomorphic to an  $i$-cell or empty.  In \cite{Che04,Che14}, we used connected and regular points to define this idea for algorithmic purposes. This is because
the concept of homeomorphism is difficult for calculation. Now we request:  that $u\cap v$ is
homeomorphic to an  $i$-cell in polynomially computable time. We also would like to restrict that idea to decide
if an $i$-cycle is a minimal cycle or an $(i+1)$-cell is also   polynomial time computable. As an example,  a polyhedron partition usually can be done in polynomial time in computational geometry.

In our definition of discrete space (a special case of one such is PL space, meaning that our definition is more strict),
a $k$-cell is a minimal closed $(k-1)$-cycle. A minimal closed $(k-1)$-cycle might not be a $k$-cell in general discrete space since it is
dependent on whether the inner part of the cell is defined in the complex or not.  We view that a $1$-cycle is a closed simple path
that is homeomorphic to a $1$-sphere. So  a $(k-1)$-cycle
is  homeomorphic to a $(k-1)$-sphere. The boundary of a $k$-cell is a $(k-1)$-cycle.


We also need another concept about regular manifolds. A regular $k$-manifold $M$ must have the following properties:
(1) Any two $k$-cells must be $(k-1)$-connected, (2) any $(k-1)$-cell must be contained in one or two
$k$-cells, (3) $M$ does not contain any $(k+1)$-cells, and (4) for any point $p$ in $M$, the neighborhood of $p$ in $M$,
denoted by $S(p)$, must be $(k-1)$-connected
in $S(M)$.

In the theory of intersection homology or PL topology~ \cite{GM},  (or as we have proved in \cite{Che13}), the neighborhood of $x$
(containing all cells that contains $x$) $S(x)$ is called the {\it star} of $x$. Note that $S(x) \setminus \{x\}$ is called the {\it link}.
Now we have: If $K$ is a piecewise linear $k$-manifold, then the link $S(x) \setminus \{x\}$ is a piecewise
linear $(k-1)$-sphere.
So we will also write ${\rm Star}(x)$ as $S(x)$ and ${\rm Link}(x)={\rm Star}(x)-\{x\}$.  In general, we can
define ${\rm Star}({\rm arc}) = \cup_{x\in {\rm arc}} {{\rm Star}(x)}$. So ${\rm Link}(\rm arc) = {\rm Star}({\rm arc}) - \{arc\}$.
${\rm Star}({\rm arc})$ is the envelope (or a type of closure) of ${\rm arc}$.

We also know that, if any $(k-1)$-cell is contained by two $k$-cells in a $k$-manifold $M$, then $M$ is closed.

In a graph, we refer to the distance as the length of the shortest path between two vertices.
The concept of {\it graph-distance} in this paper is the edge distance, meaning how many edges are needed from one vertex to another. We usually use the length of the
shortest path in between two vertices to represent the distance in graphs. In order to distinguish from the distance in Euclidean space,
we use graph-distance to represent lengths in graphs in this paper.

Therefore graph-distance is edge-distance or 1-cell-distance. It means how many 1-cells are needed to travel from $x$ to $y$.  We can generalize
this idea to define 2-cell-distance by counting how many 2-cells are needed from a point (vertex) $x$ to point $y$.  In other words,
2-cell-distance is the length of the shortest path of 2-cells that contains $x$ and $y$. In this path, each adjacent pair of
2-cells shares a 1-cell. (This path is 1-connected.)

We can define $d^{(k)}(x,y)$, the $k$-cell-distance from $x$ to $y$, as the length of the shortest path of (or the minimum of number of $k$-cells in such a sequence)
where each adjacent pair of two $k$-cells shares a $(k-1)$-cell. (This path is $(k-1)$-connected.)

We can see that $d^{(1)}(x,y)$ is the edge-distance or graph-distance. We write $d(x,y) = d^{(1)}(x,y)$

(We can also define $d^{(k)}_i(x,y)$) to be a $k$-cell path that is $i$-connected. However, we do not need to use such a concept in this paper. )}

\end{document}